\documentclass[12pt]{article}
\usepackage{bbm}
\usepackage{amsfonts}
\usepackage{graphicx} \usepackage{float}
\usepackage{amsmath}
\usepackage{mathrsfs}
\usepackage{amsthm}
\usepackage{amsmath,amssymb,amsfonts}

\renewcommand{\theequation}{\arabic{equation}}
\def\cen{\centerline}
\def\nin{\noindent}

\begin{document}

\baselineskip=20pt

\cen{{\huge { On the ruin problem in the renewal }}}
 \cen{{\huge { risk processes perturbed by diffusion}}}\vspace{1.0cm}
 \cen{{\large  Min Song}}\vspace{1.5cm}
 \cen{\large {\bf
Abstract}}
 In this paper, we consider the perturbed renewal risk process.
Systems of integro-differential equations for the Gerber-Shiu
functions at ruin caused by a claim and oscillation are established,
respectively. The explicit Laplase transforms of Gerber-Shiu
functions are obtained, while the closed form expressions for the
Gerber-Shiu functions are derived when the claim amount distribution
is from the rational family. Finally, we present numerical examples
intended to illustrate the main results.

  \vspace{0.5cm}



\vfill\noindent\footnotesize{${}^*$ Corresponding author: School of
Mathematical Sciences, Nankai University,
 Tianjin, China
\\E-mail: nksongmin@yahoo.com.cn (M. Song)

\newpage\normalsize
\baselineskip=26pt \noindent{\bf 1. Introduction} \vspace{0.4cm}
\setcounter{equation}{0}\renewcommand{\theequation}{1.\arabic{equation}}

Consider a continuous time renewal risk process perturbed by
diffusion
\begin{eqnarray}
 U(t)=u+c t -
    \sum_{i=1}^{N(t)}{ Z_{i}}+\sigma B(t), \quad t\geq0, \label{in1}
\end{eqnarray}
where $u\geq0$ is the initial capital. $c>0$ is the constant rate of
premium. The ordinary renewal process $\{N(t),t\geq0 \}$ denotes the
number of claims up to time t, with $N(t)={\text{max}}\{n\geq 1:
V_1+V_2+\cdots +V_n\leq t \}$. Then $V_i,~ i=1,2,\ldots$ are the
interclaim random variables. They are independent and assumed to
have common distribution function $K$, density function $k$, and
Laplace transforms $\hat k(s)=\int_0^{\infty}e^{-sx}k(x)dx$.
$\{Z_{i},i\geq1\}$ are independent claim-size random variables with
common distribution $P$ (such that $P(0)=0$) and density $p$.
$\{B(t),t\geq0\}$ is a standard Brownian motion with $B(0)=0$. It is
assumed that $\{N(t)\}$, $\{B(t)\}$ and $\{Z_{i}\}$ are mutually
independent and that $c E(V_i)> E(Z_i)$ providing a positive safety
loading factor.

The perturbed risk model of form (1.1) was firstly introduced by
Gerber (1970) and has been studied by many authors. See, for
example, Dufrense and Gerber (1991), Furre and Schmidli (1994),
Schmidli (1995), Gerber and Landry (1998), Wang and Wu (2000),Tsai
and Willmot (2002a, b), Zhang and Wang (2003), Li and Garrido (2005)
and the references therein.

Let $T={\text {inf}}\{t:U(t)\leq 0\}$ (with
inf$\{\emptyset\}=\infty$) be the time of ruin for risk process
(1.1). Two important nonnegative random variables in connection with
the time of ruin $T$ are $|U(T)|$, the deficit at the time of ruin,
and $U(T-)$, the surplus immediately before the time of ruin.
Consider a penalty scheme which is defined by a constant $w_0$ if
ruin occurs by oscillation and $w(U(T-),|U(T)|)$ if ruin is caused
by a jump. Then, the Gerber-Shiu discounted penalty function at ruin
$\phi(u)$ is defined by
\begin{eqnarray}
\phi(u)= \phi_w(u)+ w_0\phi_d(u), \label{in2}
\end{eqnarray}
where for $\delta\geq 0$,
\begin{eqnarray}
\phi_w(u)=E \left [e^{-\delta T}I(T<\infty,U(T)<0)w(U(T-),|U(T)|) |
U(0)=u \right ], \label{in2}
\end{eqnarray}
(with $\phi_w(0)=0$) is the expected discounted penalty function at
ruin caused by a claim, and
\begin{eqnarray}
\phi_d(u)=E[e^{-\delta T}I(T<\infty,U(T)=0)|U(0)=u], \label{in3}
\end{eqnarray}
(with $\phi_d(0)=1$) is the Laplace transform of the ruin time $T$
due to oscillation. Many ruin-related quantities can be analyzed by
appropriately choosing special penalty function $w$, for example,
let $\delta=0$ and $w=1$, then $\phi_w(u) \triangleq \psi_w(u)$
gives the probability of ruin due to a claim and $\phi_d(u)
\triangleq \psi_d(u)$ is the probability of ruin that is caused by
oscillation.

The evaluation of the Gerber-Shiu discounted penalty function, first
introduced in Gerber and Shiu(1998), is now one of the main research
problemes in ruin theory. See, for example, Gerber and Landry
(1998), Tsai and Willmot (2002a,b) for the classical surplus process
perturbed by diffusion, Li and Garrido (2005) for the generalized
Erlang (n) risk process perturbed by diffusion, Albrecher and Boxma
(2005) for the semi-Markov model, Willmot (2007) and Landriault and
Willmot (2007) for the renewal risk model, Lu and Tsai (2007) for
the Markov-Modulated process perturbed by diffusion.

The rest paper is organized as follows. In Section 2 we derive
systems of integro-differential equations for Gerber-Shiu functions.
Section 3 discusses a generalized Lundberg's equation and its roots.
And the Gerber-Shiu functions for the model are fully analyzed in
Section 4. Section 5 contains several numerical examples intended to
illustrate the main results.

 \vspace{0.5cm}
\noindent{\bf 2. Integro-differential equations} \vspace{0.4cm}
\setcounter{equation}{0}\renewcommand{\theequation}{2.\arabic{equation}}

More recently, Ren (2007) considered the risk process with
phase-type interclaim times, i.e., the distribution of the
interclaim time $K$ is phase-type with representation
$({\boldsymbol{\alpha}, \textbf{B}, \textbf{b}})$, where
$\boldsymbol{\alpha}$ and $\bf {b}$ are row vectors of length n and
$\bf {B}$ is a $n\times n$ matrix. That is, each of random variables
$V_k,k=1,2,\ldots$ corresponds to the time to absorption in a
terminating continuous-time Markov Chain $J_t^{(k)}$ with n
transient states $\{\mathcal {E}_1,\mathcal {E}_2,\ldots \mathcal
{E}_n\}$ and are absorbing state $\mathcal {E}_0$. Let $\bf e$
denote a row vector of length n with all elements being one. Then
${\bf b}^{\text T}=-{\bf B e}^{\text T}$. Following Asmussen (2000),
\begin{eqnarray*}
&&K(t)=1-\boldsymbol {\alpha} e^{{\text t}\bf B}{\bf e}^{\text T},\quad t\geq 0,\\
&&k(t)=\boldsymbol {\alpha} e^{{\text t}{\bf B}}{\bf b}^{\text
T},\quad t\geq 0,
\end{eqnarray*}
and
\begin{eqnarray}
\hat k(s)=\int_0^{\infty}e^{-st}k(t)dt=\boldsymbol {\alpha} (s\bf
I-\bf B)^{-1}{\bf b}^{\text T}.
\end{eqnarray}

For $i=1,2,\ldots,n$, let $\phi(u;i)$ denote the Gerber-Shiu
function given $U(0)=u$ and $J_0^{(1)}=\mathcal {E}_i$, that is,
\begin{eqnarray*}
\phi(u;i)=E[e^{-\delta T}I(T<\infty)w(U(T-),|U(T)|)
 | U(0)=u,J_0^{(1)}=\mathcal {E}_i],\quad i=1,2,\ldots,n.
\end{eqnarray*}
Then the Gerber-Shiu function may be computed by
\begin{eqnarray*}
\phi(u)=\boldsymbol{\alpha}~ \boldsymbol{\phi(u)}
\end{eqnarray*}
where $\boldsymbol{\phi}(u)=
(\phi(u;1),\ldots,\phi(u;n))^{\texttt{T}}$ is a column vector of
functions. Similarly, we write $\phi_w(u;i)$ and $\phi_d(u;i)$ for
the Gerber-Shiu functions at ruin caused by a claim and oscillation
respectively, given $U(0)=u$ and $J_0^{(1)}=\mathcal {E}_i$. Denoted
by $\boldsymbol{\phi}_w(u)\triangleq
(\phi_w(u;1),\ldots,\phi_w(u;n))^{\texttt{T}}$ and
$\boldsymbol{\phi}_d(u)\triangleq
(\phi_d(u;1),\ldots,\phi_d(u;n))^{\texttt{T}}$.

Our first result gives integro-differential equations for
Gerber-Shiu functions.

\noindent{\bf Theorem 2.1.}~ Let $u> 0$. Then,
$\boldsymbol{\phi}_w(u)$ satisfies the following equation
\begin{eqnarray}
\frac{\sigma ^2}{2}\boldsymbol{\phi}_w^{\prime \prime}(u)+c
\boldsymbol{\phi}_w^{\prime}(u)+ ({\bf B}-\delta~ {\bf
I})\boldsymbol{\phi}_w(u)+\left[\int_0^u \boldsymbol{\alpha}~
\boldsymbol{\phi}_w(u-x)p(x)dx+\omega(u)\right]\textbf{b}^\texttt{T}=\bf
0
\end{eqnarray}
where $\omega(u)= \int_u^{\infty}w(u,x-u)p(x)dx$, ${\bf I}=
{\text{diag}}(1,1,\ldots,1)$, $\bf 0$ denotes a column vector of
length n with all elements being $0$, and
$\boldsymbol{\phi}_w(0)={\bf 0}$, while $\boldsymbol{\phi}_d(u)$
satisfies
\begin{eqnarray}
\frac{\sigma ^2}{2}\boldsymbol{\phi}_d''(u)+c
\boldsymbol{\phi}_d'(u) -\delta \boldsymbol{\phi}_d(u)+ {\bf
B}\boldsymbol{\phi}_d(u)+\left[\int_0^u \boldsymbol{\alpha}~
\boldsymbol{\phi}_d(u-x)p(x)dx\right] \textbf{b}^\texttt{T}=\bf 0,
\end{eqnarray}
with $\boldsymbol{\phi}_d(0)=\bf e^{\text T}$.

\noindent{\bf Proof.} Define $\{J(t), t\geq 0\}$ by piecing the
$\{J_t^{(k)}\}$ together,
\begin{eqnarray*}
J(t)=\{J_t^{(1)}\},~~0\leq t<V_1, \quad
J(t)=\{J_{t-V_1}^{(2)}\},\quad V_1\leq t<V_1+V_2, \dots
\end{eqnarray*}
Then $\{J(t)\}$ is Markov. Jacobson (2005) showed that the joint
process $\{(U(t), J(t)),t\geq 0 \}$ is a Markovian additive process
and thus $\{(U(t), J(t)),t\geq 0 \}$ is a homogeneous Markov
process. Then we shall use the technique developed in Grandell
(1991, P.84) and Wu and Wei (2004) to derive the
integro-differential equations for $\phi_w(u;i)$ and $\phi_d(u;i)$
for $i=1,2,\ldots,n$. Let us consider $\phi_w(u;i)$ first. Consider
a short time interval $[0,h]$. By noting the fact that claim occurs
at the time when the state for $\{J(t), t\geq 0\}$ changes, we
separate the three possible cases:

(\emph{1}) no change of state occurs in $[0,h]$ for the underlying
Markov chain $\{J(t), t\geq 0\}$ and denoted by $J[0,h]\equiv
\mathcal {E}_i$ in this case;

(\emph{2}) the state for $\{J_t, t\geq 0\}$ changes in $[0,h]$ but
no claim occurs;

(\emph{3}) at least one claim occur in $[0,h]$,

In the case (\emph{1}). Denoted by $E^{(u,i)}$ the conditional
expectation given the initial $(U(0), J(0))=(u, \mathcal {E}_i)$ and
$\mathcal {V}(h)\triangleq u+c h+\sigma B(h)$. Let $\mathscr{F}^{J}$
and $\mathscr{F}^{(U,J)}$ denote the natural filtration of processes
$\{J(t)\}$ and $\{(U(t),J(t))\}$ respectively. For $t\geq 0$, let
$\theta _t$ be the shift operators (see, Revuz and Yor (1991,
P.34)). It follows from the Markov property of both the underlying
process $\{J(t)\}$ and the vector process $\{(U(t),J(t))\}$ that
\begin{eqnarray}
&& E^{(u,i)}\left[e^{-\delta T}I(T<\infty,U(T)<0)
w(U(T-),|U(T)|)I(J[0,h]\equiv \mathcal {E}_i)\right]\nonumber\\ &=&
E\bigg \{E\bigg [E^{(u,i)}[w (U(T-),|U(T)|)e^{-\delta
T}I(T<\infty,U(T)<0)I(J[0,h]\equiv \mathcal {E}_i) \nonumber\\
&& ~~~~~~~~
|\mathscr{F}_h^{(U,J)}]|\mathscr{F}_{\infty}^{J} \bigg ]\bigg \}\nonumber\\
&=&E \bigg\{ E\bigg[e^{-\delta h}I(J[0,h]\equiv \mathcal {E}_i)
\nonumber\\ &&  ~~~~~~~~ E^{(u,i)}\left[(w (U(T-),|U(T)|)e^{-\delta
T}I(T<\infty, U(T)<0))  \circ
\theta_h|\mathscr{F}_h^{(U,J)}\right]|\mathscr{F}_{\infty}^{J}
\bigg]\bigg\}\nonumber
\\&=&e^{-\delta h}E\bigg\{E\bigg [I(J[0,h]\equiv
\mathcal {E}_i) \nonumber\\ && ~~~~~~~~~~ E^{(U(h),J(h))}\left[w
(U(T-),|U(T)|)e^{-\delta
T}I(T<\infty,U(T)<0)\right]|\mathscr{F}_{\infty}^{J} \bigg]\bigg\}
\nonumber
\\&=&e^{-\delta h}e^{b_{ii}h}E[\phi_w(\mathcal {V}(h);i)],
\end{eqnarray}
where $b_{ij}$ is the $(i,j)$th entry of matrix $\bf B$.

For case (\emph{2}) and (\emph{3}), by the similar argument to that
of case (\emph{1}), we have
\begin{eqnarray}
&& E^{(u,i)}\left[e^{-\delta T}I(T<\infty,U(T)<0)w
(U(T-),|U(T)|)I(J(0,h]\neq \mathcal {E}_i, N(h)=0)\right]\nonumber\\
&=& e^{-\delta h}(1-e^{b_{ii}h})\sum_{j=1,j\neq
i}^{n}(\frac{b_{ij}}{-b_{ii}})E[\phi_w(\mathcal {V}(h);j)]+o(h),
\end{eqnarray}
where $o(h)/h\rightarrow 0$ as $h\rightarrow 0$, and
\begin{eqnarray}
&& E^{(u,i)}\left[e^{-\delta T}I(T<\infty,U(T)<0)w
(U(T-),|U(T)|)I(N(h)\geq 1)\right]\nonumber\\ &=& e^{-\delta
h}(1-e^{b_{ii}h})(\frac{b_i}{-b_{ii}})\sum_{j=1}^{n} \alpha_jE
\bigg[\int_0^{\mathcal {V}(h)}\phi_w(\mathcal {V}(h)-x;j)p(x)dx+
\nonumber\\&&~~~~~~~~~~~~~~~~~~~~~~~~~~~~~~~~~~~~~~   \int_{\mathcal
{V}(h)}^{\infty}w(\mathcal {V}(h),x-\mathcal {V}(h))p(x)dx
\bigg]+o(h),
\end{eqnarray}
where $b_{i}$ denotes the $i$th entry of vector $\bf b$.

Summarizing the above analysis, it follows form (2.4), (2.5) and
(2.6) that
\begin{eqnarray}
&& \phi_w(u;i)=(1-\delta h+b_{ii}h)E\left[\phi_w(\mathcal
{V}(h);i)\right]+ h\sum_{j=1,j\neq i}^{n}b_{ij}E\left[\
\phi_w(\mathcal {V}(h);j)\right]+
 h b_i \sum_{j=1}^{n}\alpha_j \nonumber\\
 &&
 E \bigg[\int_0^{\mathcal {V}(h)}\phi_w(\mathcal
{V}(h)-x;j)p(x)dx+ \int_{\mathcal {V}(h)}^{\infty}w(\mathcal
{V}(h),x-\mathcal {V}(h))p(x)dx \bigg]+o(h).
\end{eqnarray}

Apply It\^{o}'s lemma for jump-diffusion processes (McDonald (2006),
Section 20.8.) to $\phi_w(\mathcal {V}(h);k)$, for $k=1,2,\ldots, n$
we have
\begin{eqnarray*}
E[\phi_w(\mathcal {V}(h);k)]=\phi_w(u;k)+[c\phi_w
'(u;k)+\frac{\sigma ^2}{2}\phi_w''(u;k)]h+o(h).
\end{eqnarray*}

Substituting the above expressions into (2.7), canceling
$\phi_w(u;i)$ from both sides, dividing $h$ and letting
$h\rightarrow 0$ yields a system of integro-differential equations
for $\phi_w(u;i)$ given the initial surplus $u$ and the initial
state of the phase-type distribution $\mathcal {E}_i$:
\begin{eqnarray*}
&& \frac{\sigma ^2}{2}\phi_w''(u;i)+c\phi_w'(u;i)-\delta
\phi_w(u;i)+\sum_{j=1}^n b_{ij}\phi_w'(u;j)+ b_i\sum_{j=1}^n \nonumber\\
 &&~~~~~~~~~~~~~~~~~~~~~~~~~~
\alpha_j \left(\int_0^u \phi_w(u-x;j)p(x)dx+\omega
(u)\right)=0,\quad \text{for}~~ i=1,2,\ldots,n.
\end{eqnarray*}
Writing the above equations in matrix form we get (2.2) and note
that $\phi_w(0;i)=0$ for $i=1,2,\ldots,n$ since $P(T<\infty,
U(T)<0|U(0)=0)=0$. Using arguments similar to those used in deriving
(2.2), it is not difficult to get (2.3) and $\phi_d(0;i)=1$ for
$i=1,2,\ldots,n$ since $P(T<\infty, U(T)=0|U(0)=0)=1$. \quad $\Box$

\nin {\bf Remark 2.1.} When the distribution of the interclaim time
is a generalized Erlang (n) distribution,
\begin{eqnarray}
\boldsymbol{\alpha}=(1,0,\ldots,0),\bf B=\left(
\begin{array}{cccccc}
-\lambda_1 & \lambda_1  & 0         & \cdots & 0      & 0 \\
0          & -\lambda_2 & \lambda_2 &    0   &\cdots  & 0 \\
\vdots     & \vdots     & \ddots    & \ddots & \ddots & \vdots \\
0          & \ldots     &     0     &     0  &     0  & -\lambda_n
\end{array}  \right),
\bf b^{\text T}=\left(
\begin{array}{cc}
 0 \\
 0 \\
 \vdots \\
 \lambda_n
\end{array}  \right).
\end{eqnarray}
For $i=1,2,\ldots,n-1$, from (2.2) that
\begin{eqnarray*}
\lambda_i
\phi_w(u;i+1)=(\lambda_i+\delta)\phi_w(u,i)-c\phi_w'(u;i)-\frac{\sigma^2}{2}\phi_w''(u;i)
\end{eqnarray*}
and
\begin{eqnarray*}
(\lambda_n+\delta)
\phi_w(u;n)-c\phi_w'(u;n)-\frac{\sigma^2}{2}\phi_w ''(u;n)=\lambda_n
\left[\int_0^u\phi_w(u-x; 1)p(x)dx+\omega(u)\right],
\end{eqnarray*}
which are formulae (7) and (8), respectively, in Li and Garrido
(2005).

Assume that $\text{lim}_{u\rightarrow \infty}e^{-s u}{\boldsymbol
\phi}_w(u)={\bf 0}$ and $\text{lim}_{u\rightarrow \infty}e^{-s u}
{\boldsymbol \phi}_w '(u)={\bf 0}$ hold for $\Re(s)>0$. Taking
Laplace transforms on both sides of equation (2.2) and noting that
${\boldsymbol \phi}_w(0)={\bf 0}$, we have
\begin{eqnarray*}
{\bf L}(s) {\widehat {\boldsymbol{\phi}}} _w(s)=
\frac{\sigma^2}{2}{\boldsymbol{\phi}} _w '(0)- {\widehat
\omega}(s){\bf b}^{\text T},\quad s\in \mathbb{C},
\end{eqnarray*}
where
\begin{eqnarray*}
{\bf L}(s)=(\frac{\sigma^2}{2}s^2+c s-\delta){\bf I}+{\bf B}+{\bf
b}^{\text T}\boldsymbol{\alpha} \hat  p(s) ,
\end{eqnarray*}
and ${\widehat \omega}(s)=\int_0^{\infty}e^{-s u}\omega(u)du$.
Denoted by ${\bf Q}_w(s)\triangleq
\frac{\sigma^2}{2}{\boldsymbol{\phi}} _w '(0)- {\widehat
\omega}(s){\bf b}^{\text T}$. Then the vector Laplace transforms
${\widehat {\boldsymbol{\phi}}} _w(s)$ can be solved as
\begin{eqnarray*}
{\widehat {\boldsymbol{\phi}}} _w(s)=[{\bf L}(s)]^{-1} {\bf Q}_w(s).
\end{eqnarray*}
Thus
\begin{eqnarray}
{\widehat {\boldsymbol{\phi}}} _w(s)=\frac{1}{{\text {det}}[{\bf
L}(s)]}{\bf L}^{\star}(s){\bf Q}_w(s),\quad s\in \mathbb{C}.
\end{eqnarray}
where ${\bf L}^{\star}(s)$ is the adjoint of matrix ${\bf L}(s)$.

Similarly, assume that $\text{lim}_{u\rightarrow \infty}e^{-s
u}{\boldsymbol \phi}_d(u)={\bf 0}$ and $\text{lim}_{u\rightarrow
\infty}e^{-s u}{\boldsymbol \phi}_d'(u)={\bf 0}$ hold for
$\Re(s)>0$. We have
\begin{eqnarray}
{\widehat {\boldsymbol{\phi}}} _d(s)=\frac{1}{{\text {det}}[{\bf
L}(s)]}{\bf L}^{\star}(s){\bf Q}_d(s),\quad s\in \mathbb{C},
\end{eqnarray}
where ${\bf Q}_d(s)=(\frac{\sigma^2}{2}\phi _d
'(0;1)+\frac{\sigma^2}{2} s+c,\frac{\sigma^2}{2}\phi _d
'(0;2)+\frac{\sigma^2}{2} s+c,\ldots, \frac{\sigma^2}{2}\phi _d
'(0;n)+\frac{\sigma^2}{2} s+c )^{\text T}$.

It is observed from (2.9) and (2.10) that the explicit expressions
for the Gerber-Shiu functions are closely related to the roots of
equation: ${\text {det}}[{\bf L}(s)]=0$. This is discussed in the
next section.

 \vspace{0.5cm}
\noindent{\bf 3. Solutions of Lundberg's fundamental equation}
\vspace{0.4cm}
\setcounter{equation}{0}\renewcommand{\theequation}{3.\arabic{equation}}

Let $a(s)= \frac{\sigma^2}{2} s^2+cs-\delta$. For values of $s$ such
that the matrix $a(s){\bf I}+{\bf B}$ is invertible, using the same
arguments as in Ren (2007), we have
\begin{eqnarray}
{\text {det}}[{\bf L}(s)]&=&{\text {det}}[a(s){\bf I}+{\bf B}]
{\text {det}}[{\bf I}+(a(s){\bf I}+{\bf B})^{-1}{\bf b}^{\text
T}{\boldsymbol {\alpha}} \hat p(s)] \nonumber \\
 &=&{\text {det}}[a(s){\bf I}+{\bf B}]
[1+{\boldsymbol {\alpha}} (a(s){\bf I}+{\bf B})^{-1}{\bf b}^{\text
T} \hat p(s)] \nonumber \\
&=&{\text {det}}[a(s){\bf I}+{\bf B}] [1-\hat k(-a(s)) \hat p(s)],
\end{eqnarray}
where we utilize (2.1) in the last step. Since the matrix $a(s){\bf
I}+{\bf B}$ is assumed to be invertible, (3.1) indicates that the
solutions for ${\text {det}}[{\bf L}(s)]=0$ and the solutions for
Lundberg's fundamental equation
\begin{eqnarray}
\hat k(\delta-cs-\frac{\sigma^2 s^2}{2}) \hat p(s)=1,\quad s\in
\mathbb{C},
\end{eqnarray}
as defined in Gerber and Shiu (2005a) and Li and Garrido (2005) are
indentical.

\noindent{\bf Theorem 3.1.}~ For $\delta > 0$, Lundberg's
fundamental equation in (3.2) has exactly n roots, say
$\rho_1(\delta,\sigma),\rho_2(\delta,\sigma),\ldots,\rho_n(\delta,\sigma)$
with a positive real part $\Re (\rho_i(\delta,\sigma))>0$.

\noindent{\bf Proof.}~ The idea of the proof comes from Gerber and
Shiu (2005b). Let $\gamma (s)= 1/\hat k(-a(s))$. Then, as pointed
out by Gerber and Shiu (2005b), its zero occurs at $-a(s)=\xi$,
where $\xi$ ranges over all n eigenvalues of ${\bf B}$. However,
${\bf B}$ is an intensity matrix of all transient states in a
continuous time Markov chain, it is nonsingular matrix with all its
eigenvalues having negative real parts (see Corollary 8.2.1 in
Rolski et al. (1999)). Therefore, we see that $\gamma (s)$ has
exactly n positive zeros and n negative zeros.

Consider a domain that is a half disk centered at 0, lying in the
right half of the complex plane, and with a sufficiently large
radius. For $\Re(s)\geq 0$, obviously that $|\hat p(s)|\leq 1$.
Because $\gamma (s)$ has exactly n positive zeros. The theorem
follows from Rouch\'{e}'s theorem if we can show that $|\gamma
(s)|>1$ on the boundary of such a half disk.

Since phase-type distribution belongs to the rational family
distributions (see Section 4.2., below), $\gamma (s)$ has the form
of the ratios of two polynomials, in which the degree of the
denominator is less than the degree of nominator. Then $|\gamma
(s)|>1$ for $|s|$ sufficiently large. Now, for $s$ on the imaginary
axis, $\Re(s)=0$, we have $|\gamma (s)|>1/|\hat k(-a(s))|>1$ also.
This ends the proof. \quad $\Box$

\nin {\bf Remark 3.1.}  If $\delta \rightarrow 0+$ then
$\rho_i(\delta,\sigma)\rightarrow \rho_i(0,\sigma)$ for $1\leq i\leq
n$. and Eq. (3.2) becomes
\begin{eqnarray*}
\hat k(-cs-\frac{\sigma^2 s^2}{2}) \hat p(s)=1,\quad s\in
\mathbb{C}.
\end{eqnarray*}
Evidently 0 is one of the roots.

In the rest of paper, $\rho_i(\delta,\sigma)$ are simply denoted by
$\rho_i$ for $1\leq i\leq n$ and $\delta \geq 0$.

 \vspace{0.5cm}
\noindent{\bf 4. Main results } \vspace{0.4cm}
\setcounter{equation}{0}\renewcommand{\theequation}{4.\arabic{equation}}

\noindent{\textit {4.1.  Explicit expressions for $\widehat \phi$}
\vspace{0.4cm}

Here, we recall the concept of dividend differences (see Gerber and
Shiu (2005a)). For a function L(s), its dividend differences, with
respects to distinct numbers $\varrho_1,\varrho_2,\ldots$, can be
defined recursively as:
$L(s)=L(\varrho_1)+(s-\varrho_1)L[\varrho_1,s]$,
$L[\varrho_1,s]=L[\varrho_1,
\varrho_2]+(s-\varrho_2)L[\varrho_1,\varrho_2,s], \ldots$. The
definition of the dividend differences obviously can be extended to
any vector or matrix that is a function of a single variable. For
example, for matrix ${\bf L}(s)$,
\begin{eqnarray*}
{\bf L}[\varrho_1,s]&=&\frac{{\bf L}(s)-{\bf L}(\varrho_1)}{s-\varrho_1},\\
{\bf L}[\varrho_1,\varrho_2,s]&=&\frac{{\bf L}[\varrho_1,s]-{\bf
L}[\varrho_1,\varrho_2]}{s-\varrho_2},
\end{eqnarray*}
and so on. 

We assume that the roots $\rho_1,\rho_2,\ldots,\rho_n$ are distinct
in the sequel. Since $\widehat {\boldsymbol \phi}_w(s)$ is finite
for $\Re (s)\geq 0$. From (2.9) that
\begin{eqnarray*}
\frac{\sigma ^2}{2}{\bf L}^{\star }(\rho_i) {\boldsymbol
\phi}_w'(0)={\bf L}^{\star }(\rho_i){\widehat \omega}(\rho_i) {\bf
b}^{\text T},\quad {\text {for}}~ i=1,2,\ldots,n.
\end{eqnarray*}
Therefore
\begin{eqnarray*}
\frac{\sigma ^2}{2}{\bf L}^{\star }[\rho_1,\rho_2] {\boldsymbol
\phi}_w'(0)=\bigg[{\bf L}^{\star }(\rho_1){\widehat
\omega}[\rho_1,\rho_2]+{\bf L}^{\star }[\rho_1,\rho_2]{\widehat
\omega}(\rho_2)\bigg] {\bf b}^{\text T}.
\end{eqnarray*}
and recursively
\begin{eqnarray*}
\frac{\sigma ^2}{2}{\bf L}^{\star }[\rho_1,\rho_2,\ldots,\rho_n]
{\boldsymbol \phi}_w'(0)=\left[\sum_{i=1}^n{\bf L}^{\star
}[\rho_1,\rho_2,\ldots,\rho_i]{\widehat
\omega}[\rho_i,\ldots,\rho_n]\right] {\bf b}^{\text T}.
\end{eqnarray*}
Then we obtain the following results for ${\boldsymbol \phi}_w'(0)$
and ${\boldsymbol \phi}_d'(0)$:

\noindent{\bf Theorem 4.1.}  The differential of  Gerber-Shiu
functions at zero can be given by
\begin{eqnarray}
&&{\boldsymbol\phi}_w'(0) =\frac{1}{\sigma ^2/2}\{{\bf L}^{\star
}[\rho_1,\rho_2,\ldots,\rho_n]\}^{-1}\left[\sum_{i=1}^n{\bf
L}^{\star }[\rho_1,\rho_2,\ldots,\rho_i]{\widehat
\omega}[\rho_i,\ldots,\rho_n]\right] {\bf b}^{\text T}, \\&&
{\boldsymbol\phi}_d'(0) =-\frac{1}{\sigma ^2/2}\left\{{\bf L}^{\star
}[\rho_1,\rho_2,\ldots,\rho_n]\right\}^{-1}\left\{{\bf L}^{\star
}[\rho_1,\rho_2,\ldots,\rho_n]{\bf e}^{\text
T}(\frac{\sigma ^2}{2}\rho_n+c)+\right.\nonumber\\
&&\left.~~~~~~~~~~~~~~~~~~~~~~~~~~~~~~~~~~~~~~~~~~~~~~~~~~~~
 \frac{\sigma ^2}{2}{\bf L}^{\star
}[\rho_1,\rho_2,\ldots,\rho_{n-1}]{\bf e}^{\text T}\right\}.
\end{eqnarray}

Applying the dividend differences repeatedly to the numerator of Eq.
(2.9) and (2.10), we obtain the following explicit expressions for
the Laplace transform of Gerber-Shiu functions:

\noindent{\bf Theorem 4.2.} The Laplace transform of Gerber-Shiu
functions are given by
\begin{eqnarray}
&& {\widehat {\boldsymbol{\phi}}}
_w(s)=\frac{\prod_{i=1}^n(s-\rho_i)}{{\text {det}}[{\bf L}(s)]}
\bigg \{{\bf L}^{\star }[\rho_1,\rho_2,\ldots,\rho_n,s](\frac{\sigma
^2}{2}{\boldsymbol{\phi}}_w'(0)-{\widehat \omega}(\rho_n){\bf
b}^{\text T})-{\bf L}^{\star }[\rho_1, \nonumber
\\&&~~~~~~~~~~~ \ldots,\rho_{n-1},s]{\bf b}^{\text T} {\widehat \omega}[\rho_n,s]
-\sum_{i=1}^{n-1}{\bf L}^{\star }[\rho_1,\rho_2,\ldots,\rho_i]{\bf
b}^{\text T}{\widehat \omega}[\rho_i,\ldots,\rho_n,s] \bigg\},
\end{eqnarray}
and
\begin{eqnarray}
&& {\widehat {\boldsymbol{\phi}}}
_d(s)=\frac{\prod_{i=1}^n(s-\rho_i)}{{\text {det}}[{\bf
L}(s)]}\bigg\{{\bf L}^{\star
}[\rho_1,\rho_2,\ldots,\rho_n,s](\frac{\sigma
^2}{2}{\boldsymbol{\phi}}_d'(0)+(\frac{\sigma
^2}{2}\rho_n+c){\bf e}^{\text T} )\nonumber \\
&&\quad ~~~~~~~~~~~~~~~~~~~~~~~~~~~~~~~~~~~~~ +\frac{\sigma
^2}{2}{\bf L}^{\star }[\rho_1,\rho_2,\ldots,\rho_{n-1},s] {\bf
e}^{\text T}\bigg \}.
\end{eqnarray}
\noindent{\bf Proof.} By the fact that $s=\rho_1$ is a root of the
numerator in (2.9), we get
\begin{eqnarray}
{\bf L}^{\star}(s){\bf Q}_w(s)&=& {\bf L}^{\star}(s){\bf
Q}_w(s)-{\bf
L}^{\star}(\rho_1){\bf Q}(\rho_1) \nonumber \\
&=& \left[{\bf L}^{\star}(s)-{\bf L}^{\star}(\rho_1)\right]{\bf
Q}_w(s)+{\bf
L}^{\star}(\rho_1) \left[{\bf Q}_w(s)-{\bf Q}_w(\rho_1)\right]\nonumber \\
&=& (s-\rho_1)\bigg\{{\bf L}^{\star}[\rho_1,s]{\bf Q}_w(s)+{\bf
L}^{\star}(\rho_1){\bf Q}_w[\rho_1,s]\bigg\} \nonumber \\
&=& (s-\rho_1)\bigg\{{\bf L}^{\star}[\rho_1,s]{\bf Q}_w(s)-{\bf
L}^{\star}(\rho_1){\bf b}^{\text T}{\widehat
\omega}[\rho_1,s]\bigg\}.
\end{eqnarray}

Further note that $s=\rho_2$ is also a root of numerator in (2.9),
implying that $s=\rho_2$ is a zero of the expression within the
braces in (4.5), Then
\begin{eqnarray*}
{\bf L}^{\star}(s){\bf Q}_w(s)&=&(s-\rho_1){\bigg \{}\left[{\bf
L}^{\star}[\rho_1,s]{\bf Q}_w(s)-{\bf L}^{\star}(\rho_1){\bf
b}^{\text T}{\widehat \omega}[\rho_1,s]\right]\nonumber \\
&&\quad ~~~~~~~~~~~~~~~~~~~~~~~~~~~~ -\left[{\bf
L}^{\star}[\rho_1,\rho_2]{\bf Q}_w(\rho_2)-{\bf
L}^{\star}(\rho_1){\bf b}^{\text T}{\widehat
\omega}[\rho_1,\rho_2]\right]{\bigg \}}  \nonumber \\
&=&(s-\rho_1)(s-\rho_2){\bigg\{}{\bf L}^{\star}[\rho_1,\rho_2,s]{\bf
Q}_w(s)-\nonumber \\ &&\quad ~~~~~~~~~~~~~~~~~~~~~~~~~~~~~~~
\sum_{i=1}^2{\bf L}^{\star}[\rho_1,\rho_i]{\bf b}^{\text T}{\widehat
\omega}[\rho_i,\ldots,\rho_2,s]{\bigg \}},
\end{eqnarray*}
recursively from the fact $s=\rho_3,\ldots,s=\rho_{n-1}$ are roots
of the numerator in (2.9) we obtain
 \begin{eqnarray}
 {\bf L}^{\star
}(s){\bf Q}_w(s)&=&{\prod_{i=1}^{n-1}(s-\rho_i)} \bigg\{{\bf
L}^{\star }[\rho_1,\ldots,\rho_{n-1},s]{\bf Q}_w(s)\nonumber \\
&&\quad ~~~~~~~~~~~~~ -\sum_{i=1}^{n-1}{\bf L}^{\star
}[\rho_1,\rho_2,\ldots,\rho_i]{\bf b}^{\text T}{\widehat
\omega}[\rho_i,\ldots,\rho_{n-1},s] \bigg\}.
\end{eqnarray}
Further note that $s=\rho_n$ is a root of the numerator in (2.9),
from (4.6) that
\begin{eqnarray*}
{\bf L}^{\star}(s){\bf Q}_w(s)&=&{\prod_{i=1}^{n-1}(s-\rho_i)}
{\bigg \{} {\bf L}^{\star}[\rho_1,\ldots,\rho_{n-1},s]({\bf
Q}_w(s)-{\bf Q}_w(\rho_n))+ \\&& ~~~~ ({\bf L}^{\star}[\rho_1,
\ldots,\rho_{n-1},s]-{\bf
L}^{\star}[\rho_1,\ldots,\rho_{n-1},\rho_n]){\bf Q}_w(\rho_n) -\\
&&~~~~ \sum_{i=1}^{n-1}{\bf L}^{\star}[\rho_1,\ldots,\rho_i]{\bf
b}^{\text T}(\widehat {\omega}[\rho_i,\ldots,\rho_{n-1},s]-\widehat
{\omega}[\rho_i,\ldots,\rho_{n-1},\rho_n]) \bigg\}\\
&=&{\prod_{i=1}^{n-1}(s-\rho_i)}\bigg \{ {\bf
L}^{\star}[\rho_1,\ldots,\rho_{n-1},s]{\bf Q}_w[\rho_n,s]+ {\bf L}^{\star}[\rho_1, \ldots,\rho_{n},s] \\
&&~~~~~~~~~~~~~~~~~~ {\bf Q}_w(\rho_n)-\sum_{i=1}^{n-1}{\bf
L}^{\star}[\rho_1,\ldots,\rho_i]{\bf b}^{\text T} \widehat
{\omega}[\rho_i,\ldots,\rho_{n},s] \bigg\},
\end{eqnarray*}
thus formula (4.3) is derived. Formula (4.4) can be proofed in the
same way. \quad $\Box$

 \vspace{0.4cm}
 \noindent{\textit {4.2.
Closed form expressions of ${\boldsymbol{\phi}}$ for rational family
claim-size distribution}} \vspace{0.4cm}

In some cases the functions ${\boldsymbol{\phi}} _w(u)$ and
${\boldsymbol{\phi}} _d(u)$ can be explicitly and analytically
determined by inversion of (4.3) and (4.4), respectively. Consider
the case where the claim-size distribution $P$ belongs to the
rational family, i.e., its density Laplace transform is of the form
\begin{eqnarray*}
\hat p(s)=\frac{r_{m-1}(s)}{r_m(s)},\quad m\in \mathbb{N}^+,
\end{eqnarray*}
where $r_{m-1}(s)$ is a polynomial of degrees $m-1$ or less, while
$r_m(s)$ is a polynomial of degrees $m$ with only negative roots,
all have leading coefficient $1$ and satisfy $r_{m-1}(0)=r_m(0)$.
This wide class of distributions includes the Erlang, Coxian and
phase-type distributions, and also the mixtures of these (see Cohen
(1982); Tijms (1984)).

Multiply both numerator and denominator of Eq. (4.3) by $r_m(s)$,
yielding
\begin{eqnarray}
&& {\widehat {\boldsymbol{\phi}}}
_w(s)=\frac{\prod_{i=1}^n(s-\rho_i)}{r_m(s){\text {det}}[{\bf
L}(s)]}\bigg \{r_m(s){\bf L}^{\star
}[\rho_1,\rho_2,\ldots,\rho_n,s](\frac{\sigma
^2}{2}{\boldsymbol{\phi}}_w'(0)-{\widehat \omega}(\rho_n){\bf b}^{\text T})-r_m(s)\nonumber \\
 &&~~{\bf L}^{\star }[\rho_1,\ldots,\rho_{n-1},s]{\bf b}^{\text T} {\widehat \omega}[\rho_n,s] - r_m(s) \sum_{i=1}^{n-1}{\bf
L}^{\star }[\rho_1,\rho_2,\ldots,\rho_i]{\bf b}^{\text T} {\widehat
\omega}[\rho_i,\ldots,\rho_n,s] \bigg \}
\end{eqnarray}

Clearly that $r_m(s){\text {det}}[{\bf L}(s)]$ is a polynomial of
degrees $m+2n$ with leading coefficient $(\frac{\sigma ^2}{2})^n$.
So the equation $r_m(s){\text {det}}[{\bf L}(s)]=0$ has $m+2n$ roots
on the complex plane. By Theorem 3.1. and the definition of the
rational distribution, $r_m(s){\text {det}}[{\bf L}(s)]=0$ has $n$
positive real roots $\rho_1,\rho_2,\ldots \rho_n$, and $m+n$
negative real roots. Thus we can express $r_m(s){\text {det}}[{\bf
L}(s)]$ by all its roots, i.e.
\begin{eqnarray*}
r_m(s){\text {det}}[{\bf L}(s)]=(\frac{\sigma ^2}{2})^n
\prod_{i=1}^n (s-\rho_i) \prod_{i=1}^{m+n} (s+R_i),
\end{eqnarray*}
where all $R_i$' s have positive real parts. For simplicity we
assume that these  $R_i$'s are distinct. Cancel the term
$\prod_{i=1}^n (s-\rho_i)$ from both numerator and denominator of
(4.7). Consequently Eq. (4.7) can be rewritten as
\begin{eqnarray}
&& {\widehat {\boldsymbol{\phi}}} _w(s)=\frac{1}{(\frac{\sigma
^2}{2})^n \prod_{i=1}^{m+n} (s+R_i)}\bigg\{r_m(s){\bf L}^{\star
}[\rho_1,\rho_2,\ldots,\rho_n,s](\frac{\sigma
^2}{2}{\boldsymbol{\phi}}_w'(0)-{\widehat \omega}(\rho_n){\bf
b}^{\text T})-r_m(s)\nonumber \\
 &&~~~~~~~ {\bf L}^{\star
}[\rho_1,\ldots,\rho_{n-1},s]{\bf b}^{\text T}{\widehat
\omega}[\rho_n, s] - r_m(s)\sum_{i=1}^{n-1}{\bf L}^{\star
}[\rho_1,\rho_2,\ldots,\rho_i]{\bf b}^{\text T}{\widehat
\omega}[\rho_i,\ldots,\rho_n,s] \bigg\}
\end{eqnarray}

It is not difficult to see that the elements in matrix $r_m(s){\bf
L}^{\star }[\rho_1,\rho_2,\ldots,\rho_n,s]$ or $r_m(s){\bf L}^{\star
}[\rho_1,\rho_2,\ldots,\rho_{n-1},s]$ are polynomials of degrees
which are less than $m+n$. and all ${\bf L}^{\star
}[\rho_1,\rho_2,\ldots,\rho_i]$ for $i=1,2,\ldots,n$ are constants.
By decomposing the rational expressions in Eq. (4.8) into partial
fractions, we get
\begin{eqnarray*}
\frac{r_m(s){\bf L}^{\star
}[\rho_1,\rho_2,\ldots,\rho_j,s]}{\prod_{i=1}^{m+n} (s+R_i)}=
\sum_{i=1}^{m+n}\frac {{\bf M}_i^{(j )}}{s+R_i},\quad {\text {for}}~
j=n-1,n,
\end{eqnarray*}
where ${\bf M}_i^{(j )}$ for $i=1,2,\ldots,m+n$ are coefficient
matrices with
\begin{eqnarray}
{\bf M}_i^{(j )}=\frac{r_m(-R_i){\bf L}^{\star
}[\rho_1,\rho_2,\ldots,\rho_j,-R_i]}{\prod_{\ell=1,\ell\neq i}^{m+n}
(R_{\ell}-R_i)},\quad {\text {for}}~ j=n-1,n,
\end{eqnarray}
and
\begin{eqnarray*}
\frac{r_m(s)}{\prod_{i=1}^{m+n} (s+R_i)}= \sum_{i=1}^{m+n}\frac
{G_i}{s+R_i},
\end{eqnarray*}
where $G_i$ for $i=1,2,\ldots,m+n$ are coefficients given by
\begin{eqnarray}
G_i=\frac{r_m(-R_i)}{\prod_{\ell=1,\ell\neq i}^{m+n}
(R_{\ell}-R_i)}.
\end{eqnarray}

Thus, by partial fraction Eq. (4.8) can be expressed as
\begin{eqnarray}
&& {\widehat {\boldsymbol{\phi}}} _w(s)=\frac{1}{(\sigma^2/2)^n}
\sum_{i=1}^{m+n}\frac{1}{s+R_i}\bigg\{{\bf M}_i^{(n)}{\bf
Q}_w(\rho_n)-{\bf M}_i^{(n-1)}{\bf b}^{\text T}{\widehat
\omega}[\rho_n, s]\nonumber \\ &&~~~~~~~~~~~~~~~~~~~~~~~~~~~
-G_i\sum_{\ell=1}^{n-1}{\bf L}^{\star
}[\rho_1,\rho_2,\ldots,\rho_{\ell}]{\bf b}^{\text T}{\widehat
\omega}[\rho_{\ell},\ldots,\rho_n,s] \bigg\}
\end{eqnarray}
and similarly,
\begin{eqnarray}
{\widehat {\boldsymbol{\phi}}} _d(s)=\frac{1}{(\sigma^2/2)^n}
\sum_{i=1}^{m+n} \frac{1}{s+R_i} \bigg\{{\bf M}_i^{(n)}{\bf
Q}_d(\rho_n)+\frac{\sigma ^2}{2} {\bf M}_i^{(n-1)} {\bf e}^{\text T}
\bigg\}.
\end{eqnarray}

In order to determine the explicit Laplace inverse of ${\widehat
{\boldsymbol{\phi}}} _w(s)$ and ${\widehat {\boldsymbol{\phi}}}
_d(s)$, we define the same operator $T_{\mathbbm{r}}$ as in Dickson
and Hipp (2001), i.e., for an integrable real valued function $f$
with respect to a complex number $\mathbbm{r}$ ($\Re
(\mathbbm{r})\geq 0$):
\begin{eqnarray*}
T_{\mathbbm{r}}f(x)=\int_x^{\infty}e^{-\mathbbm{r}(y-x)}f(y)dy,
\quad x\geq 0.
\end{eqnarray*}
It is clear that the Laplace transform of $f$, $\hat f(s)$, can be
expressed as $T_sf(0)$, and for distinct
$\mathbbm{r}_1,\mathbbm{r}_2\in\mathbb{C}$
\begin{eqnarray*}
T_{\mathbbm{r}_1}T_{\mathbbm{r}_2}f(x)=T_{\mathbbm{r}_2}T_{\mathbbm{r}_1}f(x)=
\frac{T_{\mathbbm{r}_1}f(x)-T_{\mathbbm{r}_2}f(x)}{{\mathbbm{r}_2}-{\mathbbm{r}_1}},\quad
x\geq 0.
\end{eqnarray*}

Properties of this operator can be found in Li and Garrido (2004).
Gerber and Shiu (2005a) also presented the following useful result
on the relationship between the operator $T _{\mathbbm{r}}$ and the
corresponding dividend difference:
\begin{eqnarray}
\bigg[(\prod_{i=1}^n T_{\mathbbm{r}_i})f\bigg](0)= (-1)^{n+1} \hat
f[\mathbbm{r}_1,\mathbbm{r}_2,\ldots,\mathbbm{r}_n].
\end{eqnarray}

By (4.13) and $T_s f(0)=\hat f(s)$, Reserving the Laplace transform
of Eq. (4.11) and Eq. (4.12), we get the following theorem.

\noindent{\bf Theorem 4.3.} If the claim-size distribution belongs
to the rational family, the Gerber-Shiu functions are given by:
\begin{eqnarray}
&& {\boldsymbol{\phi}} _w(u)=\frac{1}{(\sigma^2/2)^n}
\sum_{i=1}^{m+n} \bigg\{{\bf M}_i^{(n)}{\bf Q}_w(\rho_n) e^{-R_i u}
 + e^{-R_i u} \ast \bigg[
{\bf M}_i^{(n-1)} {\bf b}^{\text T}T_{\rho_n}\omega(u) +G_i\sum_{\ell=1}^{n-1}  \nonumber \\
 && \quad ~~~~~~~~~~~~~~~~~~~~~~~~~~~~~~~ (-1)^{n-\ell} {\bf
L}^{\star }[\rho_1,\rho_2,\ldots,\rho_{\ell}]{\bf b}^{\text T}
(\prod_{j=\ell}^n T_{\rho_j}) \omega(u)   \bigg] \bigg\},
\end{eqnarray}
where $\ast$ is the convolution operator, the constants ${\bf
M}_i^{(j)}$ ($j=n-1, n$), $G_i$ and ${\bf
Q}_w(\rho_n)=\frac{\sigma^2}{2}{\boldsymbol{\phi}}_w '(0)-{\widehat
\omega}(\rho_n){\bf b}^{\text T}$ are obtained by Eq. (4.9), (4.10)
and (4.1) respectively. And
\begin{eqnarray}
{\boldsymbol{\phi}}_d(u)=\frac{1}{(\sigma^2/2)^n} \sum_{i=1}^{m+n}
e^{-R_i u} \bigg\{ {\bf M}_i^{(n)}{\bf Q}_d(\rho_n)+\frac{\sigma
^2}{2} {\bf M}_i^{(n-1)} {\bf e}^{\text T} \bigg\},
\end{eqnarray}
where ${\bf Q}_d(\rho_n)=\frac{\sigma^2}{2}{\boldsymbol{\phi}}_d
'(0)+(\frac{\sigma^2}{2} \rho_n+c){\bf e}^{\text T}$ can be
calculated from (4.2).

\vspace{0.4cm}
 \noindent{\textit {4.3.
Ruin probability}} \vspace{0.4cm}

This subsection illustrates the application of the previous results
in a special case that $\delta=0$, $w_0=1$, $w(x,y)=1$, $p(x)=\beta
e^{-\beta x}$ for $x>0$ and the interclaim times follow (2.8) with
n=2. Then (4.14) and (4.15) become the ruin probabilities caused by
a claim and oscillation respectively. We now have $\hat
p(s)=\beta/(s+\beta)$, $\omega(u)=e^{\beta u}$ for $u\geq 0$. The
matrix
\begin{eqnarray*} {\bf L}(s)=\left( \begin{array}{cc}
\frac{\sigma^2}{2} s^2+cs-\lambda_1 & \lambda_1 \\
\frac{\beta \lambda_2}{s+\beta}& \frac{\sigma^2}{2} s^2+cs-\lambda_2
\end{array}  \right)
\end{eqnarray*}
has exactly two positive real roots $\rho_1=0$, $\rho_2$ and three
negative real roots $-R_1$, $-R_2$ and $-R_3$. It follows form (4.1)
and (4.2) that
\begin{eqnarray}
{\bf Q}_w(\rho_2)&=& \frac{\sigma^2}{2}{\boldsymbol{\phi}}_w '(0)-
{\widehat \omega}(\rho_2){\bf b}^{\text T}\nonumber \\
&=& \frac{\lambda_1
\lambda_2}{\beta(\beta+\rho_2)(\frac{\sigma^2}{2}\rho_2+c)^2}
 \left(\frac{\sigma^2}{2}\rho_2+c\quad \frac{\sigma^2}{2}\rho_2+c-\frac{\lambda_2}{\beta+\rho_2}
 \right)^{\text T},\\
{\bf Q}_d(\rho_2)&=& \frac{\sigma^2}{2}{\boldsymbol{\phi}}_d '(0)+
(\frac{\sigma^2}{2}\rho_2+c){\bf b}^{\text T}\nonumber
\\&=&\frac{\lambda_1+ \lambda_2}{\rho_2+2c/\sigma^2}
 \left(1\quad  1-\frac{\lambda_2}{(\beta+\rho_2)(\frac{\sigma^2}{2}\rho_2+c)}
 \right)^{\text T}.
\end{eqnarray}
Moreover, noting that
$\psi_w(u)={\boldsymbol{\alpha}}{\boldsymbol{\psi}}_w(u)=\psi_w(u;1)$
and $\psi_d(u)=\psi_d(u;1)$, together with (4.14) and (4.15) give
the following formulae for ruin probabilities
\begin{eqnarray}
&&\psi_w(u)= \frac{1}{(\sigma^2/2)^2}\sum_{i=1}^3 \frac{\lambda_1
\lambda_2}{\beta(\beta+\rho_2){\prod_{\ell\neq i}}(R_{\ell}-R_i)}
\bigg\{(1+\frac{\frac{\sigma^2}{2}(\beta-R_i)}{\frac{\sigma^2}{2}\rho_2+c})e^{-R_i
u}-e^{-\beta u} \bigg\},\nonumber\\&&\\
&&\psi_d(u)= \sum_{i=1}^3 \frac{\beta -R_i}{{\prod_{\ell\neq
i}}(R_{\ell}-R_i)} \bigg\{(\frac{\lambda_1
+\lambda_2}{\frac{\sigma^2}{2}\rho_2+c}-R_i+\frac{2c}{\sigma^2})e^{-R_i
u} \bigg\}.
\end{eqnarray}

\vspace{0.5cm} \noindent{\bf 5. Numerical Examples} \vspace{0.4cm}
\setcounter{equation}{0}\renewcommand{\theequation}{5.\arabic{equation}}

In this section, we will present some numerical examples. In all
calculations $c=1, \sigma=1, w_0=1, p(x)=e^{-x}$ for $x\geq 0$ are
fixed. Let the interclaim times be distributed with phase-type
representation
\begin{eqnarray*}
\boldsymbol{\alpha}=(1,0),~~~~\bf B=\left(
\begin{array}{cc}
-1 & \frac{1}{2} \\
0          &-4
\end{array}  \right),~~~~
\bf b^{\text T}=\left(
\begin{array}{cc}
\frac{1}{2}\\
 4
\end{array}  \right).
\end{eqnarray*}
Then the mean $E[V_1]=9/8$, which indicates the relative safety
loading is $\frac{1}{8}$.

\noindent{\bf Example 5.1.} (Ruin probability) Let $\delta=0,
w(x,y)=1$ for $x\geq 0$, $y\geq 0$. In this case, we have
\begin{eqnarray*} {\bf L}(s)=\left( \begin{array}{cc}
\frac{s^2}{2} + s-1 +\frac{1}{2(s+1)} & \frac{1}{2} \\
\frac{4}{s+1}& \frac{s^2}{2} + s-4
\end{array}  \right),
\end{eqnarray*}
the Lundberg's equation det[{\bf L}(s)]=0 has roots:
\begin{eqnarray*}
&&\rho_1=0,~~~~\rho_2=2.06412,\\
&&{\text {and}}~~~ R_1=3.90909,~~~~R_2=3.0744,~~~~R_3=0.0806231.
\end{eqnarray*}
Then from (4.14) and (4.15) the different ruin probability
components become:
\begin{eqnarray*}
&&\psi_w(u)=0.27603 e^{-3.90909 u} -0.8912 e^{-3.0744 u} +
0.6151 e^{-0.0806231u},\\
&&\psi_d(u)=-0.2675 e^{-3.90909 u} +0.9368 e^{-3.0744 u}
+0.33066 e^{-0.0806231u},\\
&&\psi(u)=0.00853 e^{-3.90909 u} +0.0456 e^{-3.0744 u} +0.9458
e^{-0.0806231u}.
\end{eqnarray*}
\begin{figure}[H]
   \begin{center}
   \includegraphics[width=0.75\textwidth]{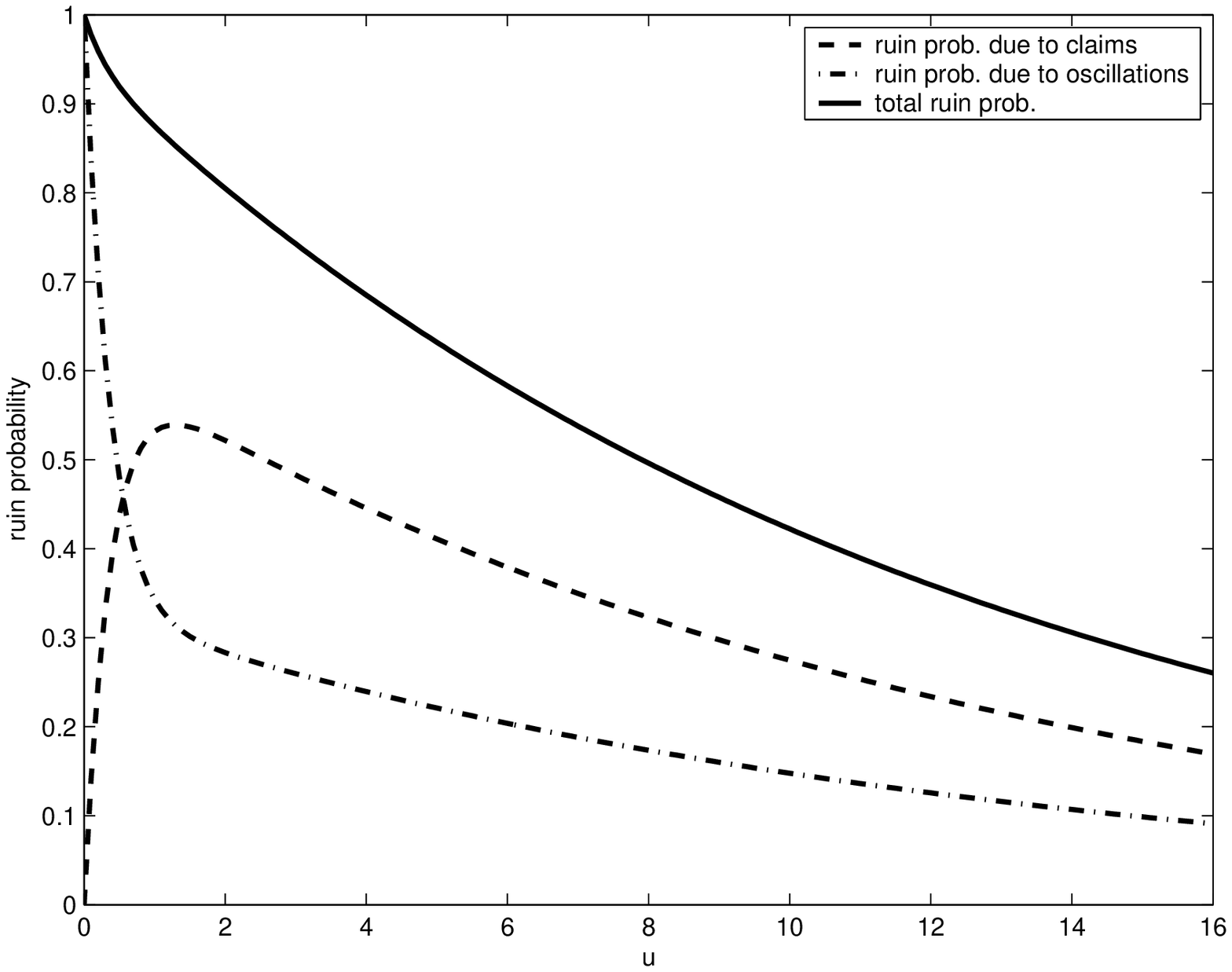}
   \cen {Fig. 1. Decomposition of the ruin probability.}
   \end{center}
   \end{figure}
With the help of Matlab, we get Fig. 1. for these ruin probabilities
for different values of u, as well as their decomposition into the
ruin probabilities due to claims and those due to oscillations. From
the graph it can be observed that the ruin probability due to
oscillations is a strictly decreasing function (from 1 to 0) of the
initial surplus u. Moreover, when u is small, it decreases sharply,
while it decreases slowly when u is large. By contrast, the ruin
probability due to claims increases quickly at first but then
deceases slowly after that.

\noindent{\bf Example 5.2.} (The Laplace transform of the ruin time)
When $w=1$,  we give Fig. 2. of $E[e^{-\delta T}I(T<\infty)|U(0)=u]$
in the case $\delta=0.1$ and $\delta=0.2$.
\begin{figure}[H]
   \begin{center}
   \includegraphics[width=0.65\textwidth]{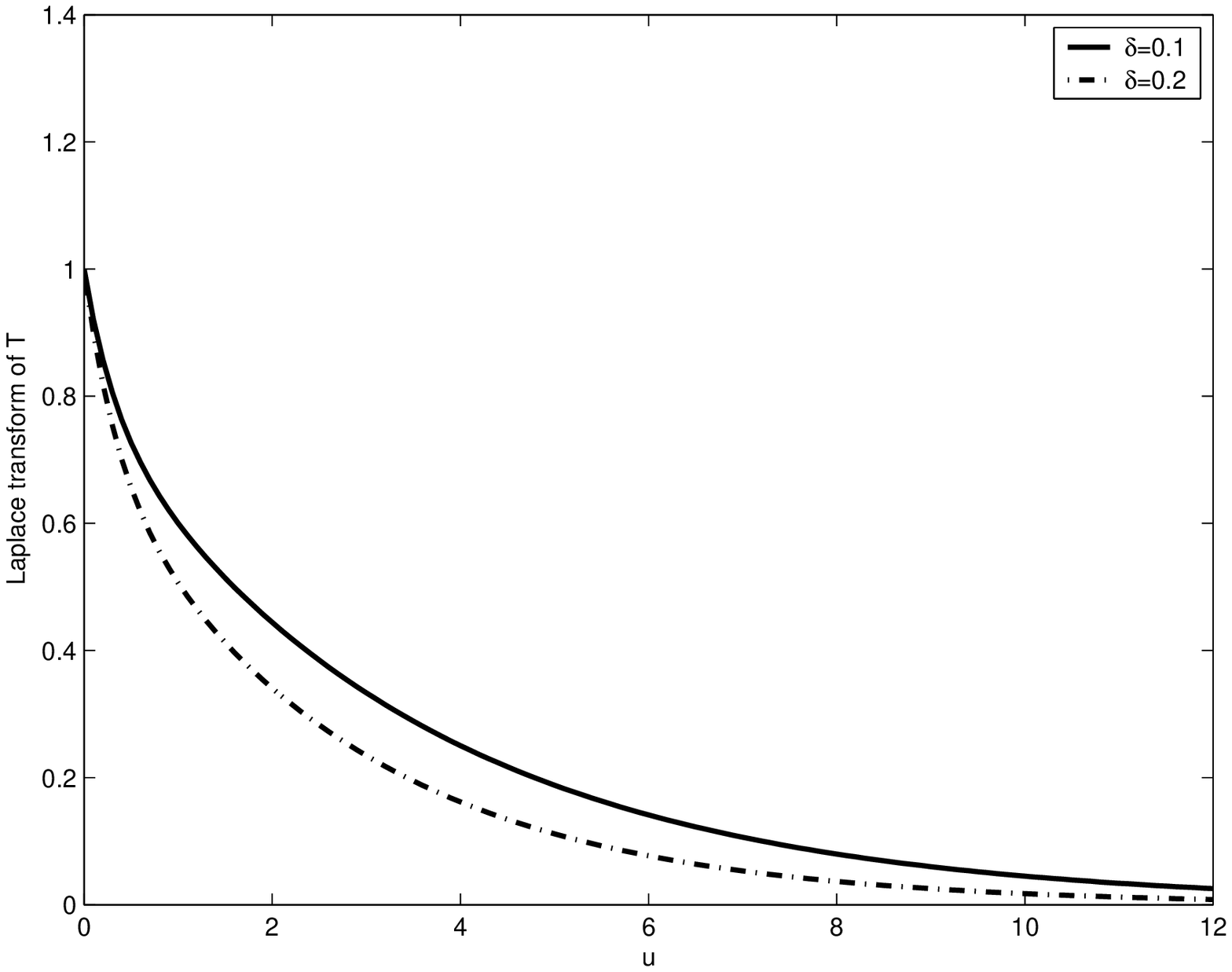}
   \cen {Fig. 2. The Laplace transform of T.}
   \end{center}
   \end{figure}
As expected, the Laplace transform of the ruin time is high for low
$\delta$, the force of interest.

\vspace{0.5 cm} \nin{\bf References} \vspace{0.4 cm}

{\small

\nin Albrecher, H., Boxma, O.J., 2005. On the discounted penalty
function in a Markov-dependent

risk model. Insurance: Mathematics and Economics 37, 650-672.

\nin Asmussen, S., 2000. Ruin Probabilities. World Scientific,
Singapore.

\nin Cohen, J.W., 1982. The Single Server Queue, Revised edition.
North-Holland, Amsterdam.

\nin Dickson, D.C.M., Hipp, C., 2001. On the time to ruin for Erlang
(2) risk process. Insurance:

Mathematics and Economics 29 (3), 333-344.

\nin Dufresne, F., Gerber, H.U., 1991. Risk theory for the compound
Poisson process that is

perturbed by diffusion. Insurance: Mathematics and Economics 10,
51-59.

\nin Furrer, H.J., Schmidli, H., 1994. Exponential inequalities for
ruin probabilities of risk pro-

cesses perturbed by diffusion. Insurance: Mathematics and Economics
15, 23-36.

\nin Gerber, H.U., 1970. An extension of the renewal equation and
its application in the collective

 theory of risk. Skandinavisk
Aktuarietidskrift, 205-210.

\nin Gerber, H.U., Landry, B., 1998. On the discounted penalty at
ruin in a jump-diffusion and

the perpetual put option. Insurance: Mathematics and Economics 22,
263-276.

\nin Gerber, H.U., Shiu,  E.S.W., 1998. On the time value of ruin.
North American Actuarial

Journal 2 (1), 48-78.

\nin Gerber, H.U., Shiu, E.S.W., 2005a. The time value of ruin in a
Sparre Andersen model.

North American Actuarial Journal 9 (2), 49-69.

\nin Gerber, H.U., Shiu, E.S.W., 2005b. Authors¡¯ Reply to
Discussions, The Time Value of Ruin

in a Sparre Andersen Model. North American Actuarial Journal 9 (2),
80-84.

\nin Grandell, J., 1991. Aspects of Risk Theory. Springer, NewYork.

\nin Jacobson, M., 2005. The time to ruin for a class of Markov
additive risk process with two-

sided jumps. Advanced in Applied Probability 37, 963-992.

\nin Landriault, D., Willmot, G.E., On the Gerber-Shiu discounted
penalty function in the Sparre

Andersen model with an arbitrary interclaim time distribution.
Insurance: Mathematics

and Economics (2007), doi:10.1016/j.insmatheco.2007.06.004

\nin Li, S., Garrido, J., 2004. On ruin for Erlang(n) risk process.
Insurance: Mathematics and

Economics 34 (3), 391-408.

\nin Li, S., Garrido, J., 2005. The Gerber-Shiu function in a Sparre
Andersen risk process per-

turbed by diffusion. Scandinavian Actuarial Journal (3), 161-186.

\nin Lu, Y., Tsai, C.C.L., 2007. The expected discounted penalty at
ruin for a Markov-Modulated

risk process perturbed by diffusion. North American Actuarial
Journal 11 (2), 136-149.

\nin McDonald, R.L. 2006. Derivatives Markets. Addison Wesley,
Boston.

\nin Ren, J.D., 2007. The discounted joint distribution of the
surplus prior to ruin and the deficit

at ruin in a Sparre Andersen model. North American Actuarial Journal
11 (3), 128-136.

\nin Revuz, D., Yor, M., 1991. Continuous Martingales and Brownian
Motion. Springer, Berlin.

\nin Rolski, T., Schmidli, H., Schmidt, V., Teugels, J.L., 1999.
Stochastic processes for insurance

and finance. John Wiley and Sons, Chichester.

\nin Schmidli, H., 1995. Cramer-Lundberg approximations for ruin
probabilities of risk processes

perturbed by diffusion. Insurance: Mathematics and Economics 16,
135-149.

\nin Tijms, H., 1994. Stochastic Models: an Algorithmic Approach.
John Wiley, Chichester.}

\nin Tsai, C.C.L., Willmot, G.E., 2002a. A generalized defective
renewal equation for the

surplus process perturbed by diffusion. Insurance: Mathematics and
Economics 30,

51-66.

\nin Tsai, C.C.L., Willmot, G.E., 2002b. On the moments of the
surplus process perturbed

by diffusion. Insurance: Mathematics and Economics 31, 327 -350.

\nin Wang, G.J., Wu, R., 2000. Some distributions for classical risk
processes that is per-

turbed by diffusion. Insurance: Mathematics and Economics 26, 15-24.

\nin Willmot, G.E., 2007. On the discounted penalty function in the
renewal risk model

with general interclaim times. Insurance: Mathematics and Economics
41 (1), 17-

31.

\nin Wu, R., Wei, L., 2004. The probability of ruin in a kind of Cox
risk model with variable

premium rate. Scandinavian Actuarial Journal (2), 121-132.

\nin Zhang, C.S., Wang, G.J., 2003. The joint density function of
three characteristics on

jump-diffusion risk process. Insurance Mathematics and Economics 32,
445-455.

 \end{document}